\documentclass[11pt, psamsfonts]{amsart}
\usepackage{amsmath}
\usepackage{amsthm}
\usepackage{amssymb}
\usepackage{amscd}
\usepackage{amsfonts}
\usepackage{amsbsy}
\usepackage{epsfig}

\theoremstyle{plain}
\newtheorem{thm}{Theorem}

\newtheorem{lem}[thm]{Lemma}

\newtheorem{prop}[thm]{Proposition}
\newtheorem{claim}[thm]{Claim}

\theoremstyle{definition}

\newtheorem{rmk}[thm]{Remark}
\newtheorem{eg}[thm]{Example}
\newcommand{\st}{such that }

\newcommand{\OOO}{\emptyset}

\newcommand{\NNN}{\mathbb{N}}

\newcommand{\QQQ}{\mathbb{Q}}
\newcommand{\RRR}{\mathbb{R}}

\newcommand{\AAA}{\mathcal{A}}
\newcommand{\BBB}{\mathcal{B}}
\newcommand{\CC}{\mathcal{C}}
\newcommand{\NN}{\mathcal{N}}
\newcommand{\MM}{\mathcal{M}}
\newcommand{\half}{\frac{1}{2}}

\begin{document}

\title[Periods, Lefschetz numbers and entropy for a class of maps
]{Periods, Lefschetz numbers and entropy for a class of maps
on a bouquet of circles}

\begin{abstract}
We consider some smooth maps on a bouquet of circles.  For these
maps we can compute the number of fixed points, the existence of
periodic points and an exact formula for topological entropy.  We
use Lefschetz fixed point theory and actions of our maps on both
the fundamental group and the first homology group.
\end{abstract}

\author[Jaume Llibre and Michael Todd]{Jaume Llibre and Michael Todd}
\address{Departament de Matem\`{a}tiques,
Universitat Aut\`{o}noma de Barcelona,\newline \indent  08193
Bellaterra, Barcelona, Spain} \email{jllibre@mat.uab.es}
\address{Mathematics Department,
University of Surrey, Guildford, \newline \indent Surrey, GU2 7XH,
UK} \email{ m.todd@surrey.ac.uk}
\thanks{The first author was partially supported by a
MCYT grant BFM2002--04236--C02--02 and by a CIRIT grant number
2001SGR 00173, the second one by a Marie Curie Fellowship number
HPMT-CT-2001-00247.}


\maketitle

\section{Introduction and statement of main results}

\label{sec:intro}

We will consider a particular class of maps on a bouquet of
circles.  We can characterise the periods of periodic orbits,
Lefschetz numbers and entropy for this class.

We first recall the concept of Lefschetz number of period $n$. Let
$M$ be a compact ANR of dimension $n$, see \cite{brown, dold}. A
continuous map $f:M \to M$ induces an endomorphism $f_{\ast k}:
H_k(M, \QQQ) \to H_k(M, \QQQ)$ for $k=0, 1, \ldots, n$ on the
rational homology of $M$. For a linear operator $A$, we let ${\rm
Tr}(A)$ denote the trace of $A$.  The {\em Lefschetz number of
$f$} is defined by
$$L(f) = \sum_{k=0}^n (-1)^k  {\rm Tr}(f_{\ast k}).$$
Since $f_{\ast k}$ are integral matrices, $L(f)$ is an integer. By
the well known Lefschetz Fixed Point Theorem, if $L(f) \neq 0$
then $f$ has a fixed point (see, for instance, \cite{brown}). We
can consider $L(f^m)$ too: $L(f^m) \neq 0$ implies that $f^m$ has
a fixed point.  However, a fixed point of $f^m$ is not necessarily
a periodic point of period $m$. Therefore, a function for
detecting the presence of periodic points of a given period was
given in \cite{l}.  This is the {\em Lefschetz number of period
$m$}, defined as
$$l(f^m)= \sum_{r|m} \mu(r) L\left(f^{\frac{m}{r}}\right),$$
where $\sum_{r|m}$ denotes the sum over all positive divisors of
$m$, and $\mu$ is the {\em Moebius function} defined as
\[
\mu(m) = \left\{\begin{array}{ll} 1 & \mbox{if $m=1$}, \\
0 & \mbox{if $k^2|m$ for some $k \in \NNN$,}\\
(-1)^r & \mbox{if $m=p_1 \cdots p_r$ for distinct prime factors.
}\end{array} \right.
\]

According to the Moebius Inversion Formula (MIF), see for example
\cite{nz},
$$L(f^m) = \sum_{r|m} l(f^r).$$
Define ${\rm Fix}(f^m)$ to be the set of fixed points of $f^m$ for
all $m \in \NNN$, and define ${\rm Per}_m (f)$ to be the set of
periodic points of period $m$.  Let ${\rm Per}(f)$ denote the set
of periods of the periodic points of $f$.

\begin{figure}[htp]
\begin{center}
\includegraphics[width=0.4\textwidth]{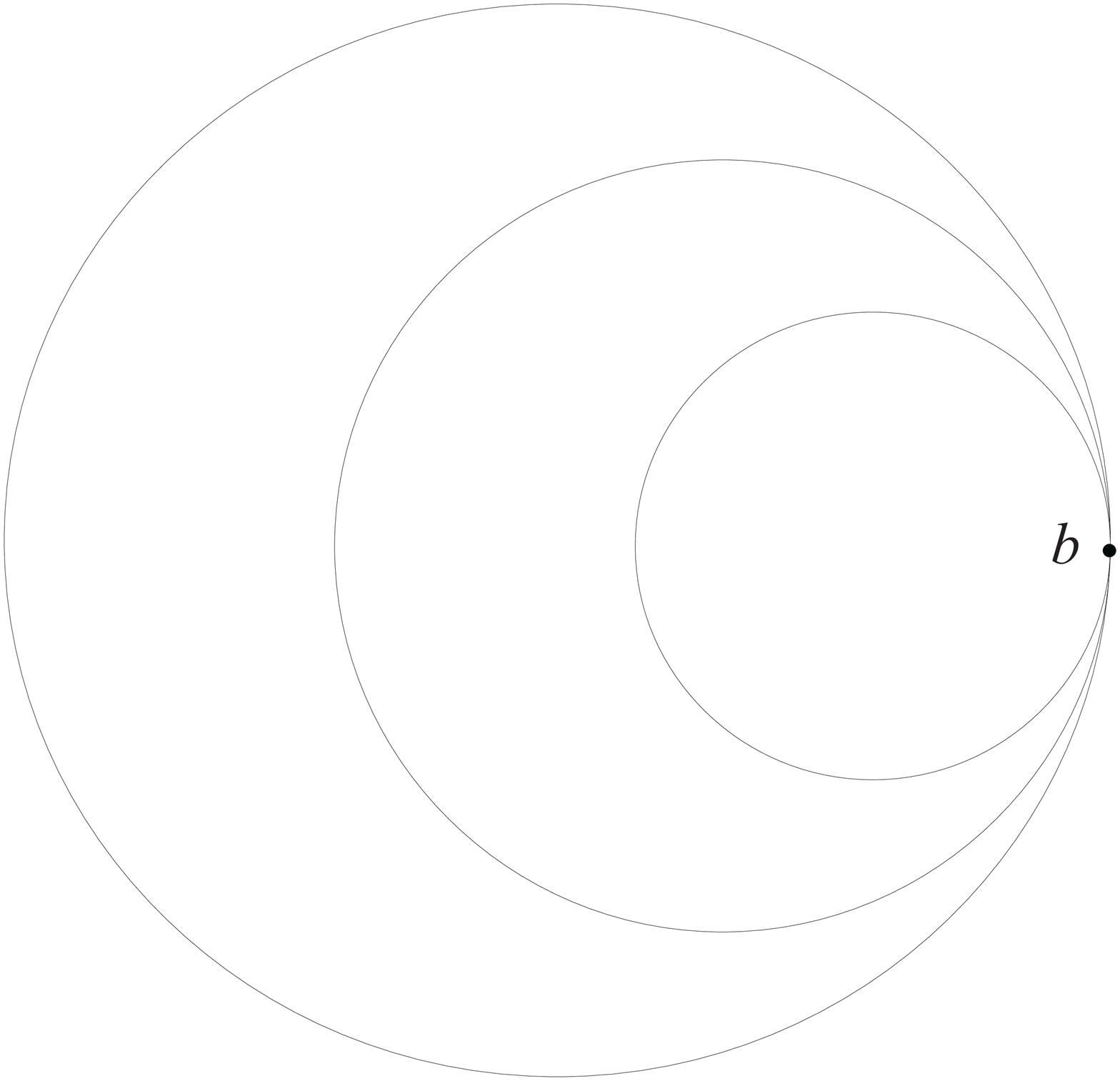} \caption{$G_3$}
\label{fig:g3}
\end{center}
\end{figure}

The Lefschetz number for periodic points has been used for
studying the set of periods of different classes of maps, see for
instance \cite{Ll}. Here we will use them for studying maps on
graphs.  In particular we consider bouquets of circles as follows.
For more details on such maps see \cite{bouquet} and
\cite{bouquet2}. We consider a set $S_1, \ldots, S_n$ in the plane
where for $1 \leq i \leq n$, $S_i$ is diffeomorphic to the unit
circle.  We call each $S_i$ a {\em circle} and suppose further
that they are nested inside each other and are pairwise disjoint,
except at a single point $b$ where they all touch. We call this
set $G_n$ and call $b$ the {\em branching point}. See
Figure~\ref{fig:g3} for a picture of some $G_3$. We give each
circle the anticlockwise orientation. With our graph arranged in
such a way, the orientation on every circle is easy to see. We say
that any graph $G \subset \RRR^n$ which is homotopic to some $G_n$
is a {\em bouquet of circles}.

Let $x, y \in S_i$ where $x \neq y$.  Then let $[x,y]$ denote the
closed arc in $S_i$ which starts at $x$, proceeds anticlockwise,
and ends at $y$. Furthermore, for $x \in S_i$ and $y \in S_j$, we
consider the connected set $[x,y]:=[x,b] \cup [b,y]$ to be an arc.
Also, we consider $\{x\}$ to be a degenerate arc.  We can extend
this definition to the open arc $(x,y)$ and the half open arcs
$(x,y]$ and $[x,y)$ in the natural way.  Note that any arc is
homotopic to a point.

Any continuous map $f:G_n \to G_n$ induces an action on $H_1(G_n,
\QQQ) = \stackrel{n}{\overbrace{\QQQ \oplus \cdots \oplus\QQQ}}$,
the first homology group.  We denote this action by $f_{\ast
1}:H_1(G_n, \QQQ) \to H_1(G_n, \QQQ)$. $f_{\ast 1}$ can be
represented by an $n \times n$ integral matrix $(m_{ij})$ \st a
generator $a_j \in H_1(G_n, \QQQ)$ maps by $f_{\ast 1}$ to the
generator $a_i$, $m_{ij}$ times, taking into account orientation.
See, for example, \cite{vick} for more details.

For a continuous map $f:M \to M$ on a compact ANR $M$, we define
the {\em minimal set of periods for $f$} to be the set $${\rm
MPer}(f)= \bigcap_{g \simeq f} {\rm Per}(g)$$ where $\simeq$
denotes homotopy.  In \cite{bouquet} and \cite{bouquet2} the
following was proved for continuous maps on $G_n$.

\begin{thm}
Let $f:G_n \to G_n$ be a continuous map and let $f_{\ast 1}$ be
the $n \times n$ integral matrix induced on the first homology
group of $G_n$.  Then the following statements hold.
\begin{itemize}
\item[(a)] If there is some element of the diagonal of $f_{\ast 1}$
different from $-2, -1,0,1,$ then ${\rm MPer}(f) = \NNN$.

\item[(b)] If all the elements of the diagonal of $f_{\ast 1}$ are $-2, -1,0$ or
$1$, and at least one of them is $-2$ then ${\rm MPer}(f) = \NNN$
or $\NNN \setminus \{2\}$.
\end{itemize}

\label{thm:lnun}

\end{thm}

Any map $f:G_n \to G_n$ has a {\em lift} to a map $\tilde{f}:[0,n]
\to [0, n]$ as follows.  We identify the integers $0, 1, \ldots,
n$ with $b$ and identify $[i-1,i)$ with $S_i$.  We assume that the
lifting map $\pi:[0,n] \to G_n$ is continuous, is orientation
preserving and is $C^1$ on each $x \in (i-1,i)$ for $1 \leq i \leq
n$. Note that $\pi$ is an example of a covering map (see
\cite{vick}).

We will consider the following class of maps, for which we can
prove more. We let $f$ be a continuous map $f:G_n \to G_n$ which
is (1) $C^1$ on $G_n \setminus \{b\}$; (2) for any $m \geq 1$, for
$x \in {\rm Fix}(f^m) \setminus \{b\}$, $|Df^m(x)|> 1$; and (3)
the sign of the derivative of the lift $\tilde{f}$, ${\rm
sign}(D\tilde{f}(x))$ for $x \in (0, n) \setminus \NNN$ is
constant.  Any such map is monotone and we say that it is in
$\MM^n$.  If, furthermore, $f^m(b) \neq b$ for all $m \geq 1$ then
we say that $f$ is in $\MM^n_b$.  Note that any $f \in \MM_b^n$ is
either, orientation preserving on all of $G_n$, or orientation
reversing on all of $G_n$.

Our first result on maps in this class is the following.

\begin{thm}
Suppose that $f \in \MM_b^n$. Then,

\begin{itemize}
\item[(a)] for all $m \geq 1$, if
$f^m$ is orientation preserving then $L(f^m)=-\# {\rm Fix}(f^m)$;
\item[(b)] for all $m \geq 1$, if
$f^m$ is orientation reversing then $L(f^m)=\# {\rm Fix}(f^m)$;
\item[(c)] if $f$ is orientation preserving then $|l(f^m)|=
\# {\rm Per}_m(f)$;
\item[(d)] if $f$ is orientation reversing and either $m$ is odd
or $4|m$, then we have $|l(f^m)|= \# {\rm Per}_m(f)$.
\end{itemize}

\label{thm:lefm}

\end{thm}

We next find formulae for the number of fixed points of maps in
terms of the action on the fundamental group.  We will find a
class of maps $\MM_\#^n$ which have an action on the fundamental
group which corresponds well with maps in $\MM^n$.

For each circle $S_j$ for $1 \leq j \leq n$ there exists a
corresponding generator in $\Pi(G_n)$. We label this generator
$a_j$.  We may assume that these are all positively oriented (that
is, each $a_j$ corresponds to a circle with anticlockwise
orientation).

We say that a word $b_1 \ldots b_m$ is {\em allowed by $\MM_\#^n$}
if either all $b_k \in \{a_1, \ldots, a_n\}$ or all $b_k \in
\{a_1^{-1}, \ldots, a_n^{-1}\}$. For a word $b_1 \ldots b_m$
allowed by $\MM_{\#}^n$, define
\[
\chi_j(b_1 \ldots b_m) = \left\{ \begin{array}{ll} \# \{b_k = a_j:
1 \leq k
\leq m \} & \mbox{if this set is not null,}\\
-\# \{b_k = a_j^{-1}: 1 \leq k \leq m \} & \mbox{if this set is
not null,}
\\
0 & \mbox{otherwise.}
\end{array}\right.
\]

Similarly we define
\[
\gamma_j(b_1 \ldots b_m) = \left\{ \begin{array}{ll} \# \{b_k =
a_j: 1 < k
< m \} & \mbox{if this set is not null,}\\
-\# \{b_k = a_j^{-1}: 1 < k < m \} & \mbox{if this set is not
null,}
\\
0 & \mbox{otherwise.}
\end{array}\right.
\]

Observe the difference between these two functions: $\chi_j$
counts the number of appearances of $a_j$ or $a_j^{-1}$ in $b_1
\ldots b_n$, but $\gamma_j$ counts the number of appearances of
$a_j$ or $a_j^{-1}$ in $b_2 \ldots b_{n-1}$. So, for example
$\chi_j(a_ja_{j+1}a_j)=2$, but $\gamma_j(a_ja_{j+1}a_j)=0$.

Now, for each $1 \leq j \leq n$, define $A_j$ to be the word
$f_\#(a_j)$.  We say that $f \in \MM_\#^n$ if all $A_j$ are
allowed by $\MM_\#^n$.  Note that $\MM^n \subset \MM_\#^n$. We
define $ d_{ij}:=\chi_i(A_j)$.

For $1 \leq k < \infty$, we say that the map $f \in \MM^n$ is in
$\MM_{b, k}^n$ if $f^k(b) = b$, but there is no $1 \leq m<k$ \st
$f^m(b) \neq b$. We say that $\MM_{b, \infty}^n= \MM_b^n$.

\begin{prop} If $f\in \MM_{b, k}^n$ for some $1 \leq k \leq
\infty$ then for any $ m  \notin k\NNN$, we have $$\#{\rm
Fix}(f^m) =
\left|1-\sum_{j=1}^n\chi_j\left(f_\#^m(a_j)\right)\right|,$$ and
if $k<\infty$, then for any $m \in k\NNN$ we have
$$\#{\rm Fix}(f^m) =
1+ \left|\sum_{j=1}^n\gamma_j\left(f_\#^m(a_j)\right)\right|.$$
\label{prop:sufficientinfo}
\end{prop}

\begin{rmk}
For our maps the action on the fundamental group and that on the
first homology group are very closely related.  However, we see by
the second part of this proposition that the fundamental group is
particularly useful when studying fixed points of maps in
$\MM_{b,k}^n$ for $k<\infty$.  In Theorem~\ref{thm:lefm} we were
not able to find an exact formula for the number of fixed points
from the Lefschetz number for maps in this class.  In fact, adding
the above result to the formula for the Lefschetz number given by
the action on the first homology group, it is possible to show
that for such maps, for $m \in k\NNN$, $L(f^m) \leq \#{\rm
Fix}(f^m) \leq 2n-1 +L(f^m)$. \label{rmk:bfixedL}
\end{rmk}

Next we prove results on periods for maps in $\MM_b^n$.

\begin{prop}
For $f\in \MM^n$, suppose that either (a) $|d_{jj}| \geq 2$ for
some $1 < j \leq n$; (b) $d_{11} \geq 2$; (c) $d_{11} < -2$; or
(d) $f \in \MM_{b, 1}^n$ and $d_{11} = -2$. Then ${\rm Per}(f) =
\NNN$. Furthermore, if (e) $d_{11}=-2$ then ${\rm Per}(f) \supset
\NNN \setminus \{2\}$. \label{prop:doubling}

\end{prop}

This is essentially the same as Theorem~\ref{thm:lnun} for maps in
$\MM^n$.  But we prove it here for completeness. We can further
characterise the set of periods in the following case.

\begin{prop}
Suppose that $f\in \MM_b^n$. Then we have the following.
\begin{itemize} \item[(a)] If there exist $1< i,j \leq n$, $i \neq j$ \st
$|d_{ij}|, |d_{ji}|\geq 1$ and $|d_{ii}|+|d_{jj}| \geq 1$, then
${\rm Per}(f) =\NNN$. \item[(b)] If there exists some $1<i \leq n$
\st $d_{i1} \neq 0, -1$, then ${\rm Per}(f) \supset \NNN \setminus
\{1\}$.
\item[(c)] If there exists some $1<i \leq n$ \st $d_{i1} = -1$,
then for all $m \geq 1$ either $m \in {\rm Per}(f)$ or $m+1 \in
{\rm Per}(f)$.
\end{itemize}
Now suppose that $f \in \MM_{b,1}^n$.  Then \begin{itemize}
\item[(d)] if there exist $1 \leq i,j \leq n$, $i \neq j$ \st
$|d_{ij}|, |d_{ji}|\geq 1$ and $|d_{ii}|+|d_{jj}| \geq 1$, then
${\rm Per}(f) =\NNN$. \end{itemize}

\label{prop:lowgrow}

\end{prop}

We next use the matrix $f_{\ast 1}$ to compute the entropy for
maps on $G_n$.  For some similar results on a different class of
maps see the recent preprint \cite{lisbon}.  We let the {\em
spectral radius} of a linear map $L$ be equal to the largest
modulus of the eigenvalues of this map. We denote this value by
$\sigma(L)$. Let $h(f)$ denote the {\em topological entropy} of
the map $f$, see Section~\ref{sec:entropy} for details.  Manning
in \cite{man} proved the following: a step towards proving the
well-known entropy conjecture, proposed by Shub in \cite{shub}.

\begin{thm}
For any continuous map $f:M \to M$, for a compact differentiable
manifold without boundary $M$, we have $h(f) \geq \log
\sigma(f_{\ast 1})$.

\label{thm:manning}
\end{thm}

Following the arguments of \cite{mislenk} we can prove the
following.  Here, given an $n \times n$ matrix $M$ we let $\|M\| :
= \sum_{i, j}|m_{ij}|$.

\begin{thm}
For a map $f\in \MM^n$ we have (a) $h(f)= \lim_{m \to \infty}
\frac{1}{m}\log \|f_{\ast 1}^m\|$; and (b) $h(f) = \log
\sigma(f_{\ast 1})$.

\label{thm:entropy}

\end{thm}

Given $f \in \MM_b^n$, the map $f_{\ast 1}$ has eigenvalues
$\lambda_1, \ldots, \lambda_d$ where the eigenvalues are in order
of decreasing modulus $|\lambda_1|\geq \cdots \geq |\lambda_d|$
(when two eigenvalues have the same modulus, any choice of order
suffices).  Our final main result is as follows.

\begin{prop}
For $f \in \MM^n$ where the eigenvalues of $f_{\ast 1}$ have
$|\lambda_1|>1$ and $|\lambda_1|> |\lambda_2|$, then there exists
some $m_0 \geq 1$ \st $m \geq m_0$ implies $m \in {\rm Per}(f)$.

\label{prop:dominant}
\end{prop}

\begin{rmk}  It should be possible to extend these results to
maps $f:G\to G$ for graphs $G$ which are homotopic to some $G_n$.
We should also be able to extend some of the results to some
classes of maps on some spaces which are homotopic to some $G_n$.
For example, some class of maps on the disk punctured $n$ times
(for maps on the twice punctured disk see \cite{eight}). However,
it is difficult to characterise such maps.
\end{rmk}

In Section~\ref{sec:lef} we prove Theorem~\ref{thm:lefm}.  In
Section~\ref{sec:fundgp} we show that the action of maps on this
class is well characterised by the action on the fundamental group
and so prove Proposition~\ref{prop:sufficientinfo}. We then go on
to prove Propositions~\ref{prop:doubling} and \ref{prop:lowgrow}.
In Section~\ref{sec:entropy} we prove Theorem~\ref{thm:entropy}.
In Section~\ref{sec:dominant} we prove
Proposition~\ref{prop:dominant}.  For examples of maps which we
can apply our results to, see Section~\ref{sec:examples}.

\section{Applying Lefschetz numbers to a bouquet of circles}

\label{sec:lef}

In this section we prove Theorem~\ref{thm:lefm} and explain the
problems associated with part (d) of the theorem.

First we recall that when $f: G_n \to G_n$ is $C^1$ and the fixed
points of $f$ are isolated, we can express
$$L(f) = \sum_{f(x)=x} {\rm ind}(f, x),$$ where ${\rm ind}(f, x)$ is
the {\em index of $f$ at $x$}.   If $x \neq b$ then ${\rm ind}(f,
x) = (-1)^{u_+(x)}$, where $u_+(x)=1$ whenever $Df(x)>1$ and
$u_+(x)=0$ otherwise.  For more details see \cite{kelly} or
\cite{lnun}. There, the question of the index of $f$ at $b$ when
$b$ is a fixed point is also discussed.

\begin{proof}[Proof of Theorem~\ref{thm:lefm}] The first two statements
of the theorem are easy to see because, $$L(f) = \sum_{f(x)=x}
{\rm ind}(f, x)= \sum_{f(x)=x} (-1)^{u_+(x)}$$ where $u_+(x)$ is
defined as above.

So $L(f)$ counts the number of fixed points, giving negative or
positive sign if $f$ is orientation preserving or reversing,
respectively.  So we have proved (a) and (b).

Next we prove (c).  Since $f$ is orientation preserving, the
summands for $l(f^m)$ are all negative. Therefore, by the MIF,
$$\sum_{r|m} \# {\rm Per}_{r}(f) = |L(f^m)|.$$
From the definition of $l(f)$, applying the MIF again we have
$|l(f^m)|= \# {\rm Per}_{m}(f)$.

To prove (d), we first suppose that $m$ is odd.  Then the summands
for $l(f^m)$ are of the form $\mu(r) L(f^{\frac{m}{r}})$ where
$r|m$. Since $\frac{m}{r}$ cannot be even, $L(f^{\frac{m}{r}})$
are either all negative or all positive depending on whether $f$
is orientation preserving or reversing, respectively. Therefore,
by the MIF,
$$\sum_{r|m} \# {\rm Per}_{r}(f) = |L(f^m)|.$$
Again, a further application of the MIF gives $|l(f^m)|= \# {\rm
Per}_{m}(f)$.

Now if $4 |m$ then let $n$ be \st $m=4n$.  Any summand for
$l(f^m)$ is of one of the following forms.
\begin{itemize}
\item[(i)] $\mu(r) L(f^{\frac{4n}{r}})$ where $r|4n$ and $r$ is odd
(so $r|n$).  Since $\frac{4n}{r}$ is even, $L(f^{\frac{4n}{r}})$
is negative in any case.

\item[(ii)] $\mu(2r)L(f^{\frac{2n}{r}})$ where $r|2n$ and $r$ is odd
(so $r|n$). Since $\frac{2n}{r}$ is even, $L(f^{\frac{2n}{r}})$ is
negative in any case.

\item[(iii)] $\mu(4r)L(f^{\frac{n}{r}})$ where $r|n$. Since $\mu(4r)=0$,
this term is null.

\end{itemize}

Thus all of the terms $L(f^{\frac{m}{r}})$ which contribute to
$l(f^m)$ are negative and so, applying the MIF as above we see
that
$$|l(f^m)|= \# {\rm Per}_{m}(f).$$ \end{proof}

\begin{rmk}
We cannot extend this method directly to maps with attracting
periodic points, even if they are monotone.  For example, we can
create a monotone $C^1$ map which has every repelling fixed point
followed by an attracting one.  So we can have $L(f)=0$, where $f$
has arbitrarily many fixed points. (If $f$ is orientation
reversing this does not make any difference for $L(f)$. But this
presents a problem for $L(f^2)$.)
\end{rmk}

\begin{rmk}
We explain why this result cannot be extended to $m$ where $2 |m$,
but $4\nmid m$ when $f$ is orientation reversing.  Suppose that
$m=2p$ for some $p$ prime (we obtain similar problems if $m = 2p_1
\ldots p_r$ with $p_i >2$ prime). Since $f$ is orientation
reversing,
\begin{eqnarray*} l(f^{2p}) & = & L(f^{2p}) -
L(f^{p})- L(f^{2}) + L(f) \\
& = & -\left[\#{\rm Per}_{2p}(f) + \#{\rm Per}_{p}(f)+ \#{\rm
Per}_{2}(f) + \#{\rm Per}_{1}(f)\right] \\ & & - \left[\#{\rm
Per}_{p}(f) +\# {\rm Per}_{1}(f)\right]  + \left[\#{\rm
Per}_{2}(f) +
\#{\rm Per}_{1}(f)\right] + \#{\rm Per}_{1}(f) \\
& = & -\#{\rm Per}_{2p}(f)- 2\#{\rm Per}_{p}(f).
\end{eqnarray*}
So, even when $l(f^m) \neq 0$, we cannot be so sure about the
presence of periodic points of period $m$.  This is seen in the
following examples.

\label{rmk:even}

\end{rmk}

It is convenient to construct our examples on the level of
homology where we only have information about $f_{\ast 1}:H_1(G_n,
\QQQ) \to H_1(G_n, \QQQ)$.

\begin{figure}[htp]
\begin{center}
\includegraphics[width=0.5\textwidth]{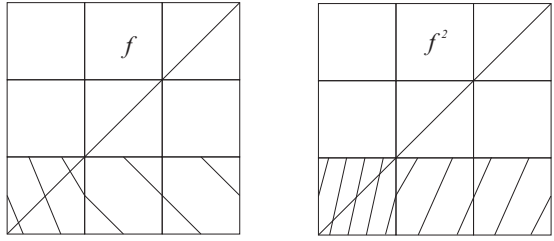}
\caption{Lift for a map $f \in \MM_b^3$ where $l(f^2) \neq 0$, but
$2 \not\in {\rm Per}(f)$.} \label{fig:even}
\end{center}
\end{figure}

\begin{eg}
Consider the map $f \in \MM^1$ which has action $f_{\ast 1}$ on
$H_1(G_1, \QQQ)$ equal to multiplication by $m_{11}$ where
$m_{11}=-2$. Then $f_{\ast 1}^2$ is multiplication by $4$. We
calculate $L(f)= 3$, $L(f^2) =-3$. So $l(f^2) =-6$.  By
Remark~\ref{rmk:even} we have
$$l(f^2) = -\#{\rm Per}_{2}(f)- 2\#{\rm Per}_{1}(f).$$ So we deduce
that ${\rm Per}_{2}(f) = \OOO$.  Therefore, the Lefschetz number
for periodic points does not always detect periodic points of even
order when the original map is orientation reversing.

We can also construct further such examples for any $n \geq 2$ as
follows.  See Figure~\ref{fig:even} for an example on $G_3$. Let
$(m_{ij})$ be the matrix representing the action of $f_{\ast 1}$
on $H_1(G_n, \QQQ)$.  Now suppose that $m_{11} =-2$; $m_{1j} = -1$
for $1 \leq j \leq n$; and $m_{ij} =0$ for $0<i\leq n, 1 \leq j
\leq n$.  Here we obtain the same behaviour on $S_1$ as on $G_1$
above (note that there are no periodic points outside $S_1$ here).
Therefore, we cannot be sure in such cases that $l(f^m) \neq 0$
implies that there are periodic points of period $m$.

\end{eg}

\section{Finding periods from the action on the fundamental group}

\label{sec:fundgp}

In fact, most of the information on periodic points for maps in
$\MM^n$ can be read from the action on the fundamental group
$\Pi(G_n)$. We will see that there is a one to one correspondence
between maps with a particular type of action on $\Pi(G_n)$ and
homology classes of maps in $\MM^n$. (As we will note later, this
is not the case when we consider the action on first homology.)

\subsection{Coding of $f$ on the fundamental group}

\label{ssec:coding}



If a word $b_1 \ldots b_m$ is allowed by $\MM_\#^n$ and has
$(\{b_1\} \cup \{b_m\}) \cap \left(\{a_1\} \cup
\{a_1^{-1}\}\right) \neq \OOO$ then we say that $b_1 \ldots b_m$
{\em is allowed by $\MM_{\# b}^n$}.  Note that a map with the
action $f_{\#} : a_{j} \mapsto a_{j_1}\ldots a_{j_{n_j}}$ for
$j_{k} \in \{1, \ldots, n \}$ which starts and finishes at the
same point in the circle corresponding to $a_{j_1}$, is homotopic
to a map with the action $f_{\#} : a_j \mapsto a_{j_2}\ldots
a_{j_{n_j}}a_{j_1}$. We can argue analogously in the orientation
reversing case.  So for maps in $\MM_\#^n$ with $f^k(b) \neq b$
for all $k \geq1$, we may assume that $a_{j_1} \in \left\{a_1,
a_1^{-1}\right\}$ for all $1 \leq j \leq n$. Observe that the
action $a_j \mapsto a_{1} a_{j_2} \ldots a_{j_{n_j}}$ gives an
orientation preserving map on $G_n$ which starts at $f(b)$, then
covers the arc $[f(b), b]$; then covers in turn the circles
$S_{j_2}, \ldots, S_{j_{n_j-1}}$ and $S_{j_{n_j}}$; finally it
covers the arc $[b, f(b)]$.

\begin{lem}
Suppose that $A_1, \ldots A_n$ are allowed by $\MM_{\# b}^n$ and
if $d_{j1} =0$ for all $1<j \leq n$ then $d_{11} \neq -1$. Then
there exists a map $f \in \MM^n$ with the action $f_\#:a_j \to
A_j$ for $1 \leq j \leq n$.  Furthermore, any $g \in \MM_\#^n$
with the same action is homotopic to $f$. \label{lem:fundgp}
\end{lem}

\begin{proof}
We will find a piecewise linear lift map $\tilde{g}: [0, n] \to
[0,n]$ with the required action and then show that $\tilde{f}:[0,
n] \to [0,n]$, the lift of $f$ must be homotopic to $\tilde{g}$.

For an interval $J$, and a linear map $g:J \to \RRR$, let $|Dg|_J
= |Dg(x)|$ for any $x \in J$. Given $1 \leq i \leq n$, we consider
the word $A_i$.  We let $\tilde{g}:[0,n] \to [0,n]$ be the
piecewise linear map with $\tilde{g}(j) = \half$ for $0 \leq j
\leq n$; $|D\tilde{g}|_{(j-1, j)} = {\rm sign}(\chi_1(a_1))n_j$;
first the map has $\tilde{g}(j-1) = \half$; then it covers half of
$[0,1]$ before covering the intervals $[i-1, i]$ given in $A_j$ in
the order given by $A_j$; finally the map covers
$\left[0,\half\right]$ where $\tilde{g}(j) = \half$.

For example if $f: G_3 \to G_3$ and $f_{\#}: a_1 \mapsto
a_1a_3a_1a_2a_2$ then $D\tilde{g}|_{[0,1]} = 5$ and has
$\tilde{g}\left(\left[0, \frac{1}{10}\right)\right) =
\left[\frac{1}{2}, 1\right)$, $\tilde{g}\left(\left[\frac{1}{10},
\frac{3}{10} \right)\right) = \left[2, 3\right)$,
$\tilde{g}\left(\left[\frac{3}{10}, \frac{5}{10} \right)\right) =
\left[0, 1\right)$, $\tilde{g}\left(\left[\frac{5}{10},
\frac{7}{10} \right)\right) = \left[1, 2\right)$,
$\tilde{g}\left(\left[\frac{7}{10}, \frac{9}{10} \right)\right) =
\left[1, 2\right)$, and $\tilde{g}\left(\left[\frac{9}{10}, 1
\right)\right) = \left[0, \frac{1}{2}\right)$.

We now show that $\tilde{g}$ is homotopic to $\tilde{f}$.  We
first may assume that $\tilde{f}$  has been `pulled  tight'.  That
is, we choose a homotopy  which results in a local homeomorphism,
i.e. given any $0 \leq j \leq n$,  for all $x \in (j-1,j)$ there
exists a neighbourhood $U$ of $x$ \st $\tilde{f}|_U$ is a
homeomorphism. This means that the graph of $f$ has no null
homotopic loops.

Suppose that $f$ is orientation preserving.  For $1 \leq j \leq
n$, let $I_{j1}$ be the minimal interval in $[0,1]$ \st
$\tilde{f}:I_{j1} \to [\tilde{f}(0), 1]$ is a surjection.  Let
$\hat{I}_{j1}$ be the equivalent interval for $\tilde{g}$. Since
$\tilde{f}$ is assumed to be a local homeomorphism,
$\tilde{f}_{I_{j1}}$ is a homeomorphism.  Since
$\tilde{f}_{I_{j1}}$ and $\tilde{g}_{\hat{I}_{j1}}$ are both
homeomorphisms on intervals with the same orientation then they
are homotopic.

Now for any small enough interval $U$ adjacent and to the right of
$I_{j1}$ we claim that $\tilde{f}(U) \subset [j_2-1, j_2]$.  If
not then there is some $i \neq j_2$ \st $\tilde{f}(U) \subset
[i-1, i]$.  But since $\tilde{f}$ is a local homeomorphism, we can
extend $U$ so that $\tilde{f}(U) = [i-1, i]$.  But then $a_{j_2} =
a_i$ which is a contradiction. As above, we can show that
$\tilde{f}_{I_{j2}}$ and $\tilde{g}_{\hat{I}_{j2}}$ are homotopic.
We may continue this process up to $n_j$ to prove that $\tilde{f}$
and $\tilde{g}$ are homotopic.

Next we need to show that $\tilde{g}$ gives a map $g:G_n \to G_n$
which is in $\MM^n$.  We need to show that for any fixed point $x$
of $\tilde{g}^m$, we have $|D\tilde{g}^m(x)|>1$. We fix some $1
\leq j \leq n$ and consider $(j-1,j)$. We have the following
cases.

{\bf Case 1:} There exists some $i \neq 1$ \st $|d_{ij}| \geq 1$.
Then $|D\tilde{g}|_{(j-1, j)} \geq 2$.

{\bf Case 2:} Suppose that we are not in Case 1.

{\bf Case 2a:} Suppose $j \neq 1$. Then since we are not in Case
1, there are no fixed points of $\tilde{g}$ in $[j, j-1]$. The
only way to obtain fixed points is to take some iterate
$\tilde{g}^m$ which passes through some interval $[i-1, i]$ which
has $[j-1,j]$ in its image. The interval $(i-1,i)$ must be in Case
1, so we have $|D\tilde{g}|_{[i-1,1]} \geq 2$. Therefore,
$|D\tilde{g}^m|_{(j-1, j) } \geq 2$.

{\bf Case 2b:} Suppose $j=1$.  If $d_{11}=1$ and $d_{1i} =0$ for
all $1 < i \leq n$ then we proceed similarly to Case 2a since we
do not have any fixed points in $[0,1]$.

If $d_{11} =-1$ and $d_{i1} \leq 1$ for some $1 \leq i \leq n$
then again we have $|D \tilde{g}|_{[0,1]} \geq 2$.

Therefore, in all cases for $x \in {\rm Fix}(\tilde{g}^m)$,
$|D\tilde{g}^m(x)|>1$.

Letting $g := \pi \tilde{g} \pi^{-1}$, we are finished.
\end{proof}

Note that we can often find some homotopic map $f$ which is also
in $\MM_b^n$.

\begin{rmk}
Suppose that $f \in \MM_{\# b}^n$, $f^k(b) \neq b$ for all $k \geq
1$, $d_{11}=-1$ and, contrary to Lemma~\ref{lem:fundgp},
$d_{j1}=0$ for all $1 < j \leq n$.  Then $f$ has two fixed points
$x_1, x_2$ in $S_1$.  It is easy see that a piecewise linear
version on $f|_{S_1}$ would have $|Df(x_1)|, |Df(x_2)| = 1$, so
this map could not be in $\MM^n$. It is possible in some cases to
perturb so that $|Df(x_1)|, |Df(x_2)|
>1$, but this will always create some points $y \in S_1$ with
$|Df(y)|<1$ which could mean that there are points $x \in G_n$
with $f^m(x) = x$ and $|Df^m(x)| <1$, i.e. $f \notin \MM^n$.
\label{rmk:degen}

\end{rmk}

\begin{proof}[Proof of Proposition~\ref{prop:sufficientinfo}] We first
suppose that $f \in \MM_b^n$.  We consider the lift $\tilde{f}$.
If $f$ is orientation preserving then
\begin{equation} \#{\rm Fix}(f)= -1+ \sum_{j=1}^n \chi_j(A_j).
\label{eqn:chi} \end{equation} (We will explain this in our case,
but it can also be seen for $G_2$ by looking at the proof of
Proposition 2 of \cite{lnun}). The reason for this is that for
$1<j \leq n$, the image of $\tilde{f}([j-1,j])$ will start at
$\tilde{f}(b)$ and, if there is some $i$ \st $a_{j_k} = a_j$, then
this image must start below the diagonal $\{(x, x): 0 \leq x \leq
n\}$ and cross it in order to cover $[j-1,j]$.  This gives a fixed
point every time this crossing happens.

When we are dealing with $\tilde{f}$ on $[0,1]$ we note that our
map must miss the diagonal following the first appearance of $a_1$
in $A_1$.  But for every subsequent appearance of $a_1$ there is a
corresponding fixed point (the $-1$ term in (\ref{eqn:chi})
accounts for this).

If $f$ is orientation reversing then
\begin{equation} \#{\rm Fix}(f)= 1-\sum_{j=1}^n \chi_j(A_j).
\label{eqn:chi2}\end{equation} This is essentially the same as the
orientation preserving case except that any image
$\tilde{f}([0,1])$ must cross the diagonal as many times as $a_1$
occurs in $A_1$, {\em plus 1}.  See for example
Figure~\ref{fig:even}.

Clearly, given any $m \geq 1$, we can replace $A_j$ with
$f_\#^m(a_j)$ in (\ref{eqn:chi}) or (\ref{eqn:chi2}) to find
$\#{\rm Fix}(f^m)$ as required.

For $f \in \MM_{b,k}^n$ where $k< \infty$ and $m \in \NNN
\setminus k\NNN$ then the proof is the same as above.  Now suppose
that $m \in k\NNN$.  For any $1 \leq j \leq n$, if $f_\#^m(a_j)$
has first or last element equal to $a_j$ then the graph of
$\tilde{f}$ on $(j-1,j)$ has no corresponding crossing of the
diagonal. However, if $a_j$ appears anywhere else in $f_\#^m(a_j)$
there is a corresponding crossing of the diagonal. Hence there are
$|\gamma_j(f_\#^m(a_j))|$ fixed points of $f^m$ in $S_j \setminus
\{b\}$.  By assumption, there is also a fixed point of $f^m$ at
$b$, so
$$\#{\rm Fix}(f^m)= 1+ \sum_{j=1}^n |\gamma_j(f_\#^m(a_j))|$$ as
required. \end{proof}

By Proposition~\ref{prop:sufficientinfo} the set of fixed points
of $f\in \MM^n$ are completely determined by the action on the
fundamental group.

\subsection{Finding periodic points from the fundamental group
action} \label{ssec:perfund}

\begin{proof}[Proof of Proposition~\ref{prop:doubling}] In all cases,
$|d_{jj}| \geq 2$ for some $1 \leq j \leq n$, and so $f$ has a
fixed point in $S_j$. We will show that in cases (a), (b), (c) and
(d), when $m \geq 1$, $\#{\rm Fix}(f^{m+1}|_{S_j})$, the number of
fixed points of $f^{m+1}$ in $S_j$, is greater than $\#{\rm
Fix}(f^{m}|_{S_j})$, the number of fixed points of $f^m$.
Therefore, there must be some new fixed point of $f^{m+1}$, which
has not been counted before as a fixed point for any $f^p$ where
$p \leq m$.  Hence we must have a periodic point of period $m+1$
in $S_j$.  Since this will be true for any $m \geq 1$, we have
${\rm Per}(f) = \NNN$. In case (e) this argument will follow for
any $m \geq 2$ and so ${\rm Per}(f) \supset \NNN \setminus \{2\}$.

{\bf Case 1:}  First suppose that  $f \in \MM_{b, k}^n$ where
$k=1$, i.e. $f(b)=b$ and we are in case (d).  Then we can see from
the proof of Proposition~\ref{prop:sufficientinfo} that the number
of fixed points of $f^p$ in $S_j$ is $1+|\gamma_j(f_\#^p(a_j))|$
(note that the 1 counts the fixed point at $b$).  Therefore if we
can show that
\begin{equation} 1+|\gamma_j(f_\#^{m+1}(a_j))| >
1+|\gamma_j(f_\#^{m}(a_j))| \label{eqn:fixedgrowth} \end{equation}
then $\#{\rm Fix}(f^{m+1}|_{S_j})> \#{\rm Fix}(f^{m}|_{S_j})$ and
there must be a periodic point of period $m+1$ in $S_j$.

Every element $a_j$ in the word $f_\#^p(a_j)$ gives rise to two
occurrences of $a_j$ in $f_\#^{p+1}(a_j)$.  Therefore,
$|\gamma_j(f_\#^{p+1}(a_j))| \geq 2|\gamma_j(f_\#^{p}(a_j))|$, so
(\ref{eqn:fixedgrowth}) is satisfied whenever
$|\gamma_j(f_\#^{p}(a_j))|>1$.  Since $|d_{jj}| \geq 2$ this is
true for any $p>1$, so ${\rm Per}(f) \supset \NNN \setminus
\{2\}$. For $p=1$ we have three cases: (i) if
$\gamma_j(f_\#(a_j))=0$ then $|\gamma_j(f_\#^2(a_j))| \geq 2$, so
$\#{\rm Fix}(f^2|_{S_j})> \#{\rm Fix}(f|_{S_j})$; (ii) if
$|\gamma_j(f_\#(a_j))|=1$ then $|\gamma_j(f_\#^2(a_j))| \geq 3$,
so $\#{\rm Fix}(f^2|_{S_j})> \#{\rm Fix}(f|_{S_j})$; (iii) if
$|\gamma_j(f_\#(a_j))|=2$ then $|\gamma_j(f_\#^2(a_j))| \geq 4$,
so $\#{\rm Fix}(f^2|_{S_j})> \#{\rm Fix}(f|_{S_j})$.  Therefore,
in all these cases, ${\rm Per}(f) = \NNN$.

From now on we will assume that $k>1$.

{\bf Case 2:} We consider $m<k-1$.  The proof here also follows
when $f \in \MM_b^n$.  We can see from the proof of
Proposition~\ref{prop:sufficientinfo} that for $p \leq m$, $$ -1+
|\chi_j(f_\#^p(a_j))| \leq \#{\rm Fix}(f^p|_{S_j}) \leq 1+
|\chi_j(f_\#^p(a_j))|.$$  Since $|d_{jj}| \geq 2$,
$|\chi_j(f_\#^{p+1}(a_j))| \geq 2|\chi_j(f_\#^{p}(a_j))|$ for any
$p \geq 1$.  Hence, we have
$$\#{\rm Fix}(f^{m+1}|_{S_j}) \geq -1+ 2|\chi_j(f_\#^{m}(a_j))|
\geq \#{\rm Fix}(f^{m}|_{S_j}) + |\chi_j(f_\#^{m}(a_j))|-2.$$
Therefore, $\#{\rm Fix}(f^{m+1}|_{S_j})> \#{\rm
Fix}(f^{m}|_{S_j})$ whenever $|\chi_j(f_\#^{m}(a_j))|>2$.  This is
always the case for $m \geq 2$.

For $m = 1$ we have \begin{equation} \#{\rm Fix}(f^2|_{S_j}) \geq
-1+ 2|\chi_j(f_\#(a_j))|. \label{eqn:lowterm} \end{equation} If we
are in case (a) then we have $\#{\rm Fix}(f|_{S_j})=
|\chi_j(f_\#(a_j))|$, so by (\ref{eqn:lowterm}), $\#{\rm
Fix}(f^2|_{S_j})
> \#{\rm Fix}(f|_{S_j})$ and we are finished. For
cases (b), (c) and (e), we have $j=1$. If we are in case (b) then
$\#{\rm Fix}(f|_{S_1})= |\chi_1(f_\#(a_1))|-1$, so
 by (\ref{eqn:lowterm}), $\#{\rm Fix}(f^2|_{S_1})
> \#{\rm Fix}(f|_{S_1})$ and we are finished.  If we are in case
(c) then $\#{\rm Fix}(f|_{S_1})= |\chi_1(f_\#(a_1))|+1$ and
$|\chi_1(f_\#(a_1))|>2$, so  by (\ref{eqn:lowterm}), $\#{\rm
Fix}(f^2|_{S_1}) > \#{\rm Fix}(f|_{S_1})$ and we are finished.
Note that case (d) is covered in Case 1.  In case (e) it is
possible that $\#{\rm Fix}(f^2|_{S_1}) = \#{\rm Fix}(f|_{S_1})$,
so we can't be sure if $2$ is a period or not.

This proof also follows for the rest of the set $\NNN \setminus
\{p: p =lk \ {\rm or} \ p = lk -1 \ {\rm for } \ l \in \NNN\}$.

{\bf Case 3:} We consider $m=lk-1$ for any $l \in \NNN$ where
$k<\infty$. Considering the formulas for the number of fixed
points of $f^p$ in $S_j$, it is sufficient to show that
\begin{equation} 1+|\gamma_j(f_\#^{m+1}(a_j))| >
|\chi_j(f_\#^{m}(a_j))|+1 \label{eqn:growthatk}. \end{equation} We
 compute that for a word $b_1 \ldots b_n$ allowed by $\MM_\#^n$,
$\gamma_j(b_1 \ldots b_n) \geq |\chi_j(b_1 \ldots b_n)| -2$.
Therefore,
$$|\gamma_j(f_\#^{m+1}(a_j))| \geq 2|\chi_j(f_\#^{m}(a_j))|-2.$$
Hence, if $|\chi_j(f_\#^{m}(a_j))| > 2$ then (\ref{eqn:growthatk})
is satisfied.  Since $|\chi_j(f_\#^{m}(a_j))| \geq 2^{m}$,
(\ref{eqn:growthatk}) holds for case (c), or whenever $m>1$.  For
$m=1$, we compute that in case (a), $j=1$ and
\begin{eqnarray*} \# {\rm Fix}(f^2|_{S_1}) & = &
1+|\gamma_1(f_\#^{m+1}(a_1))| \geq
-1+|\chi_1(f_\#^2(a_1))| \\
& \geq &  -1+2|\chi_1(f_\#(a_1))| = \# {\rm Fix}(f|_{S_1}) -1+
|\chi_1(f_\#(a_1))|.\end{eqnarray*}  Since $|\chi_1(f_\#(a_1))|
\geq 2$ we are finished.  In case (b), $j=1$
\begin{eqnarray*} \# {\rm Fix}(f^2|_{S_1}) & = &
1+|\gamma_1(f_\#^{m+1}(a_1))| \geq -1+|\chi_1(f_\#^2(a_1))| \\
& \geq & -1+2|\chi_1(f_\#(a_1))| = \# {\rm Fix}(f|_{S_1})+
|\chi_1(f_\#(a_1))|.\end{eqnarray*}  Since $|\chi_j(f_\#(a_j))|
\geq 2$ we are finished.

{\bf Case 4:} We consider $m = lk$ for any $l \in \NNN$ where
$1<k<\infty$.  Similarly to above, it is sufficient to show that
\begin{equation} |\chi_j(f_\#^{m+1}(a_j))|-1 >
|\gamma_j(f_\#^{m}(a_j))|+1 \label{eqn:growthafterk}.
\end{equation}

We have $$|\chi_j(f_\#^{m+1}(a_j))| \geq 2|\chi_j(f_\#^{m}(a_j))|
\geq |\chi_j(f_\#^{m}(a_j))|+ |\gamma_j(f_\#^{m}(a_j))|.$$
Therefore, (\ref{eqn:growthafterk}) is satisfied whenever
$|\chi_j(f_\#^{m}(a_j))| >2$.  But this is always true when $k>1$.
\end{proof}

\begin{rmk}
As mentioned above, any map $f \in \MM^n$ has a matrix action
$(m_{ij})$ on the first homology which has either all entries
positive or all entries negative. Here the terms $m_{ij}$ take the
place of $d_{ij}$. This was considered in \cite{bouquet2} and with
$b$ fixed in \cite{bouquet}. There the proof of the final part of
Theorem~\ref{thm:lefm} was proved applying Bolzano's Theorem to
subgraphs.

Note that it is not the case that any map with such an action is
homotopic to a map in $\MM^n$.  This is because the action on the
first homology abelianises the action on the fundamental group.
So, in particular, there exist homotopy classes with this action
on the first homology, for which every map in the class has
positive local degree at some point and negative local degree at
some other point. (We say a map has {\em positive (negative) local
degree} if the map is locally orientation preserving (reversing).)
\end{rmk}

\begin{proof}[Proof of Proposition~\ref{prop:lowgrow}] Suppose first that
we are in case (a). We suppose that $|d_{jj}| \geq 1$. Since we
also have $|d_{ij}|, |d_{ji}| \geq 1$, we have
$|\chi_j(f_\#^m(a_j))| \geq 1$ for all $m \geq 0$ and
$|\chi_j(f_\#^m(a_i))| \geq 1$ for all $m \geq 1$.  Note that in
particular, $f$ has a fixed point in $S_j$.

Let $A_j \mapsto a_1 a_ja_i$. (This is the simplest case for
$j>1$, and, in terms of creating periodic points, the worst since
$a_j$ only appears once in $A_j$.) Then we prove that every
application of $f_{\#}$ to $f_{\#}^m(a_j)$ creates a {\em new}
fixed point in $S_j$.  We will use the fact that $f_{\#}$ is a
homomorphism repeatedly.

We have $f_{\#}^{m+1}(a_j) = f_{\#}^m(a_1) f_{\#}^m(a_i)
f_{\#}^m(a_j)$. The function $\chi_j$ counts the number of times
$a_j$ appears in a given word.  Thus, $|\chi_j(f_\#^{m+1}(a_j))| =
|\chi_j(f_{\#}^m(a_1) f_{\#}^m(a_i) f_{\#}^m(a_j))| \geq
|\chi_j(f_\#^{m}(a_j))|+1$ since $|\chi_j(f_\#^m(a_i))| \geq 1$
for all $m \geq 0$. Therefore, the application of $f_\#$ generates
a new fixed point in $S_j$. Whence $m \in {\rm Per}(f)$ for all $m
\geq 1$. Note that if $A_j$ is a longer word, we obtain the same
result (in that case, the number of fixed points created by each
iteration could be even greater). Furthermore, if $f$ is
orientation reversing we can apply the same proof.  The proof of
(d) follows similarly.

If we are in case (b) and not in case (a) or a case covered by
Proposition~\ref{prop:doubling} then we are in the orientation
preserving case. The simplest form for $A_1$ is $a_1 a_i$. We have
$f_\#^{m+1}(a_1) = f_{\#}^m(a_1) f_{\#}^m(a_i)$. So again the fact
that $|\chi_1(f_\#^m(a_i))| \geq 1$ for all $m \geq 1$ means that
we have found a new fixed point in $S_1$. Whence $m \in {\rm
Per}(f)$ for all $m > 1$.

If we are in case (c), then we find our new fixed points in $S_i$.
We may suppose that $A_1= a_1^{-1} a_i^{-1}$ and $A_i=a_1^{-1}$.
Then
$$f_{\#}^{m+1}(a_i) = f_{\#}^{m}\left(a_1^{-1}\right) =
f_{\#}^{m-1}\left(a_1a_i\right) =f_{\#}^{m-1}\left(a_1\right)
f_{\#}^{m-1}\left(a_i\right).$$ Since
$\left|\chi_i\left(f_{\#}^{m}\left(a_1\right)\right)\right| \geq
1$ for $m \geq 1$, we find a new fixed point in $S_i$ after the
application of $f^2$.  Furthermore, each subsequent application of
$f^2$ yields a new fixed point. \end{proof}

See Example~\ref{eg:lowgrow} for an application of this. The
following is an easy corollary of Proposition~\ref{prop:lowgrow}.

\begin{prop}
Suppose that there exists some $m> 1$ \st for $f \in \MM_\#^n$,
replacing $d_{ij}$ with $\chi_i(f_\#^m(a_j))$ and $\MM_{b,1}^n$ by
$\MM_{b,m}^n$ in (d), we satisfy the conditions of
Proposition~\ref{prop:lowgrow}. Then we have the same conclusions
when we replace $\NNN$ with $m\NNN$ (in (b), the conclusion is
replaced by ${\rm Per}(f) \supset m(\NNN \setminus \{1\})$).
\label{prop:delaylowgrow}

\end{prop}

See Example~\ref{eg:delaylowgrow} for an application of this.

\begin{rmk}
Note that given a map  $f \in \MM_b^n$, this map has the minimal
number of periodic points within the class of maps which are
homotopic to $f$ and which have $b$ non--periodic (it is shown in
\cite{bouquet} that the maps which minimise the number of fixed
points within this homotopy class have $b$ fixed). For example, if
$a_1 \mapsto a_1a_3a_1a_2a_2$, then $f$ must cross the diagonal
the number of times $a_1$ appears in the action minus 1 (minus one
because $f(b) \in S_1$ implies that we go from $f(b)$ to $f(b)$ in
$S_1$ without crossing the diagonal exactly once). But this is
precisely what our maps do, and no more.  We can argue similarly
for $f \in \MM^n_{b, k}$.

\end{rmk}

\section{Periodic points and entropy}

\label{sec:entropy}

The main aim of this section is to prove
Theorem~\ref{thm:entropy}. This involves showing that the
eigenvalues of the matrices $f_{\ast 1}^m$ give us a lot of
information about periodic points and about entropy. We first give
entropy in terms of a limit involving $f_{\ast 1}^m$, proving
Theorem~\ref{thm:entropy}(a) and then, for part (b), we give
entropy in terms of the spectral radius of $f_{\ast 1}$.

We will give some basic definitions for entropy, see, for example
\cite{jaumebook} for more details. Let $X$ be a compact Hausdorff
metric space. We say that the set $\AAA$ is an {\em open cover}
for $X$ if $\bigcup_{A \in \AAA} A \supset X$ and all $A$ are open
sets. A {\em subcover} of $X$ from $\AAA$ is a subset of $\AAA$
which is also a cover of $X$.  When it is clear what $X$ is, we
simply refer to covers and subcovers.

Let $\AAA$ be an open cover of $X$ . For a continuous map $f:X \to
X$, we define $f^{-i}(\AAA) := \{f^{-i}(A): A \in \AAA\}$,
$\bigvee_{i=1}^{m-1} \AAA_i := \{A_1 \cap \ldots \cap A_{m-1}: A_i
\in \AAA_i, {\rm and }\ A_1 \cap \ldots \cap A_{m-1} \neq \OOO\}$
and $\AAA^m := \bigvee_{i=0}^{m-1} f^{-i}(\AAA)$. Also let
$\NN(\AAA)$ be the minimal cardinality of any subcover from
$\AAA$.

Let $$h(f, \AAA) := \lim_{m \to \infty} \frac{1}{m} \log
\NN(A^m).$$  Then we define the {\em topological entropy of $f$}
to be
$$h(f) := \sup h(f, \AAA)$$ where the supremum is taken over all
open covers of $X$.

Here we will let $X$ be some bouquet $G_n$. We say that $\AAA$ is
a {\em cover of $G_n$ by arcs} if $\bigcup_{A\in \AAA } A \supset
G_n$, each $A \in \AAA$ is an arc of $G_n$ and all $A \in \AAA$
are pairwise disjoint. These arcs can be open or closed or half
open and half closed or even degenerate.  (This notion is similar
to `a cover by intervals' when the phase space is the interval,
see Chapter 4.2 of \cite{jaumebook}.)

Let $\AAA$ be a cover of $G_n$. We call $\AAA$ an {\em $f$--mono
cover} if for all $A_i \in \AAA$ there is some circle $S_j$ \st
$f:A_i \to S_j$ is an injective homeomorphism. Note that if $\AAA$
is an $f$--mono cover then $\AAA^m$ is an $f^m$--mono cover.

The following results will allow us to prove
Theorem~\ref{thm:entropy}(a). Propositions~\ref{prop:arccover} and
\ref{prop:entropy} are adapted versions of the theory of
\cite{mislenk}.  We follow the exposition of this theory in
\cite{jaumebook}.

\begin{prop}
For $f \in \MM_n$, $h(f) = \sup(f, \AAA)$ where the supremum is
taken over finite covers by arcs.

\label{prop:arccover}
\end{prop}

For the proof of this see Proposition 4.2.2 of \cite{jaumebook}
which proves that this is so for interval maps and finite covers
by intervals.

\begin{lem}
Suppose that $f \in \MM^n$.  Suppose that $M$ is the matrix
$f_{\ast 1}$.  Then there is a natural $f$--mono cover by arcs
with cardinality $\|A\|$.

\label{lem:mincover}

\end{lem}

\begin{proof}
We construct the cover as follows.  Considering the lift
$\tilde{f}:[0,n] \to [0,n]$, let $P = \tilde{f}^{-1}(b) = \{p_1,
\ldots p_m\}$  where $p_1< \cdots < p_m$.  For $1 \leq i <m$, let
$P_i = [p_i, p_{i+1})$.  Also, let $P_m = [p_m, n] \cup [0, p_1)$.
Let $A_i = \pi(P_i)$ (note that $\pi(P_m)$ is an arc) and let
$\AAA=\{A_1, \ldots, A_m\}$.

Since $\#(\AAA) = \#(\tilde{f}^{-1}(b))$, we have $\#(\AAA) =
\|A\|$ as required.
\end{proof}

The following is Proposition 4.2.3 of \cite{jaumebook} with minor
adaptations so that it applies to our case (we must adapt the
situation for interval maps to the situation for maps on
bouquets). We include a proof for completeness.

\begin{prop}
For $f\in \MM^n$, if $\AAA$ is a mono--cover of $G_n$ then
$h(f)=h(f, \AAA)$.

\label{prop:entropy}

\end{prop}

\begin{proof} Let $\tilde{\BBB}$ be a finite cover of $G_n$ by
arcs.  Let $\BBB = \tilde{\BBB} \vee \AAA$.  Let $\CC$ be a cover
chosen from $\AAA^m$.  Take $A \in \CC$.  The map $f^k|_{A}$ is a
homeomorphism for $k= 1, \ldots , m$.  Therefore, for any $B \in
\BBB$, the set $A \cap f^{-k}(B)$ is an arc (unless it is empty).
Each arc has at most $2$ endpoints (note that a degenerate arc has
only one endpoint). Let $x$ be an endpoint of an element $B \in
\BBB^m$. Then there exists $A \in \CC$ \st $A \cap B \neq \OOO$
and $x$ is a endpoint of $A \cap B$.  Since $B =
\bigcap_{k=0}^{m-1}f^{-k}(B_k)$ for some $B_0, \ldots, B_{m-1} \in
\BBB$ and each of the sets $f^{-k}(B_k)$ is a union of a finite
number of arcs, $x$ is an endpoint of some component of
$f^{-k}(B_k)$ for some $k \in \{0, \ldots, m-1\}$.  Hence $x$ is
an endpoint of $A \cap f^{-k}(B)$ for this $k$.  In each $A\in
\CC$ there are at most $2m\# (\BBB)$ such endpoints.  The number
of possible arcs with endpoints in a given set is not larger than
$4$ times the square of the cardinality of this set (we multiply
by 4 because arcs with given endpoints may or may not contain
them). Therefore, $\#( \BBB^m|_{A}) \leq 4(2m\# (\BBB^m))^2 $.
Hence, $\NN(\BBB^m) \leq 4(2m\# (\BBB^m))^2 \# (\CC)$. Since $\CC$
was arbitrary, we obtain
$$\NN(\BBB^m) \leq 4(2m\# (\BBB^m))^2 \NN(\AAA^m).$$ In the
limit we get $$h(f, \tilde{\BBB}) \leq h(f, \BBB) \leq h(f,
\AAA).$$

By Proposition~\ref{prop:arccover}, in calculating the entropy we
need only consider finite covers by arcs, so we have $h(f) \leq
h(f, \AAA)$, and consequently $h(f) = h(f, \AAA)$.
\end{proof}

The following is proved in the appendix of \cite{milngen}.

\begin{lem} For a matrix of complex numbers $M$, the limit
$\lim_{m \to \infty}\log\|M^m\|^{\frac{1}{m}}$ exists.

\label{lem:limitmat}

\end{lem}

\begin{proof}[Proof of Theorem~\ref{thm:entropy}(a)] Consider the
$f$--mono cover $\AAA$ of $G_n$ constructed in
Lemma~\ref{lem:mincover}.  We let $M$ be the action of $f_{\ast
1}$ on the first homology.  Then $\# (\AAA) = \|M\|$. Furthermore,
$\# (\AAA^m )= \|M^m\|$. Therefore, by
Proposition~\ref{prop:entropy},
$$h(f) = h(f, \AAA) = \lim_{m \to \infty}\frac{1}{m} \log\|M^m\|.$$
Since, by Lemma~\ref{lem:limitmat} this limit exists (we could
also refer to Section 4.1 of \cite{jaumebook} to show that any
such limit of the cardinality of the pullback of covers exists),
Theorem~\ref{thm:entropy}(a) is proved. \end{proof}

The proof of Theorem~\ref{thm:entropy}(b) is a simple corollary of
Theorem~\ref{thm:entropy}(a) and the following result: Theorem A.3
of \cite{milngen}.  The proof also follows from \cite{lisbon}.

\begin{thm}
The spectral radius of any real or complex $n \times n$ matrix $M$
is given by
$$\sigma(M) = \lim_{k \to \infty}\|M^k\|^{\frac{1}{k}} =
\lim\sup_{k \to \infty}\|{\rm Tr}(M^k)\|^{\frac{1}{k}}.$$

\label{thm:spmattr}

\end{thm}

For applications of Theorem~\ref{thm:entropy}, see any of the
examples in Section~\ref{sec:examples}.

\section{Computing periods from eigenvalues of $f_{\ast 1}$}
\label{sec:dominant}

As above, the spectral radius of $f_{\ast 1}$ can be computed as
$\lim\sup_{k \to \infty}\|{\rm Tr}(M^k)\|^{\frac{1}{k}}$. However,
if $\lim_{k \to \infty}\|{\rm Tr}(M^k)\|^{\frac{1}{k}}$ exists
then we can say more.  We first state a result of \cite{asymplef}
(in fact, there the theorem also extends to maps with higher
homologies than we consider here).  We need the following
definition. A $C^1$ map $f:M \to M$ of a compact $C^1$
differentiable manifold is called {\rm transversal} if $f(M)
\subset M$ and for all $m \in \NNN$, for all $x \in {\rm
Fix}(f^m)$ we have det$(I-df^m(x)) \neq 0$, i.e. 1 is not an
eigenvalue of $df^m(x)$.

\begin{thm}
Let $M$ be a compact manifold with $H_i(M, \QQQ) =0$ for $i>1$.
Suppose that $f:M \to {\rm Int}(M)$ is a $C^1$ transversal map.
Further, assume that the limits
$$\lim_{m \to \infty}\left|{\rm Tr}\left(f_{\ast 1}^m\right)
\right|^{\frac{1}{m}}$$ and
$$\lim_{m \to \infty}\left|\sum_{d|m}\mu(d){\rm Tr}
\left(f_{\ast 1}\right)\right|^{\frac{1}{m}}$$ exist.  If there is
an eigenvalue different from a root of unity or zero then there
exists $m_0\geq 1$ \st
\begin{itemize}
\item[(a)] for all $m\geq m_0$ odd we have that $m \in {\rm
Per}(f)$;
\item[(b)] for all $m\geq m_0$ even we have that $\{\frac{m}{2},m\}
\cap {\rm Per}(f) \neq \OOO$.

\end{itemize}

\label{thm:asymplef}
\end{thm}

\begin{rmk}
Suppose that $|\lambda_1|>1$ and $|\lambda_1|>|\lambda_2|$.Then we
claim that the limit $\lim_{m \to \infty}\left|{\rm
Tr}\left(f_{\ast 1}^m\right) \right|^{\frac{1}{m}}$ exists and is
equal to $|\lambda_1|$ since
$$|\lambda_1|^k-(d-1)|\lambda_2|^k<|{\rm Tr}\left(f_{\ast 1}^k\right)|<
|\lambda_1|^k+(d-1)|\lambda_2|^k.$$ Taking limits we prove the
claim.

In a similar way we are able to show that
$$\lim_{m \to \infty}\left|\sum_{d|m}\mu(d){\rm Tr}
\left(f_{\ast 1}\right)\right|^{\frac{1}{m}}$$ exists.  This is
because it can be shown that there exists some $C>0$ \st
$$|\lambda_1|^k - \frac{1}{C}m|\lambda_1|^{\frac{k}{2}}<
\left|\sum_{d|m}\mu(d){\rm Tr} \left(f_{\ast 1}\right)\right| <
|\lambda_1|^k +Cm|\lambda_1|^{\frac{k}{2}}.$$ (For a calculation
of this type, see the proof of Proposition~\ref{prop:dominant}.)
Letting $k \to \infty$ we again obtain $|\lambda_1|$ as the limit.
Therefore we have the conclusions of Theorem~\ref{thm:asymplef}
for our map.

\label{rmk:limtr}

\end{rmk}

In fact, in our class, we can improve this result to obtain
Proposition~\ref{prop:dominant}.  For our proof, we need to show
that if a particular growth condition on the number of fixed
points is satisfied then we can be sure of the existence of some
periodic points.  To give an idea of this approach we state the
following easily proved claim.

\begin{claim}
Suppose that  $f:M \to M$ is some map on some space $M$.  If we
have $$\#{\rm Fix}_m(f) > \sum_{r|m, r \neq m}\#{\rm Per}_r(f)$$
then $$\# {\rm Per}(f^m) = \#{\rm Fix}_m(f) - \sum_{r|m, r \neq
m}\#{\rm Per}_r(f)>0.$$ \label{claim:growth}

\end{claim}

Now we give the main tool for the proof of
Proposition~\ref{prop:dominant}.

\begin{prop}
Suppose that  $f:M \to M$ is some map on some space $M$.  If for
some $m \geq 1$,
$$\# {\rm Fix}(f^m)> \sum_{\frac{m}{k} {\rm prime}, k \neq m} \# {\rm
Fix}(f^k)$$ then ${\rm Per}_m(f) \neq \OOO$. \label{prop:fmbig}
\end{prop}

\begin{proof}
We have $$\#{\rm Fix}(f^m) = \sum_{r|m} \# {\rm Per}_r(f).$$  Now
supposing that ${\rm Per}_m(f) = \OOO$,

$$\# {\rm Fix}(f^m) = \sum_{r|m, r \neq m} \# {\rm Per}_r(f).$$

So if we prove that
$$\sum_{\frac{m}{k} {\rm prime}, k \neq m} \# {\rm Fix}(f^k) \geq
\sum_{r|m, r \neq m} \# {\rm Per}_r(f),$$ then the proposition
will follow.

Note that we can write $m$ as a product of prime factors $m = p_1
\ldots p_n$.  Thus, if $\frac{m}{k}$ is  prime and $k \neq m$ then
$k=\frac{p_1 \ldots p_n}{p_i}$ for some $1 \leq i \leq n$. So
$$
\# {\rm Fix}(f^k)   =   \# {\rm Fix}\left(f^{\frac{p_1 \ldots
p_n}{p_i}}\right)  = \sum \#{\rm Per}_{r}(f)$$

where the sum runs over all combinations $r= p_{j_1} \cdots
p_{j_{n_r}}$ where $j_1 < \cdots <j_q$ and all $j_k \in \{1,
\ldots, n\} \setminus \{i\}$.

We can express any $r|m$ which has $r \neq m$ as prime factors:
$r= p_{j_1} \ldots p_{j_{n_r}}$ where $1 \leq n_r<n$.  Therefore,
the term
$\# {\rm Per}_r (f)$ is counted $\left(\begin{array}{cc} n-1  \\
 n_r
\end{array}\right)$ $(\geq 1)$ times by the sum on the left, but only once by
the sum on the right.  So the proposition is proved.
\end{proof}

\begin{proof}[Proof of Proposition~\ref{prop:dominant}.] First suppose
that $f \in \MM_b^n$.  We assume that for some $m> 1$, ${\rm
Per}_m(f) = \OOO$. Otherwise we are finished. From
Proposition~\ref{prop:fmbig}, to prove that we contradict this
assumption on periodic points, it is sufficient to show that
\begin{equation} \# {\rm Fix}(f^m)> \sum_{\frac{m}{k} {\rm prime}, k \neq m}
\# {\rm Fix}(f^k).\label{eqn:fmgrow} \end{equation}

By Theorem~\ref{thm:lefm}, for any $m \geq 1$, $\#{\rm
Fix}(f^m)=|L(f^m)| = |1-(\lambda_1^m + \cdots + \lambda_d^m)|$.
Clearly, for $m \geq 2$,
$$m\left(|\lambda_1|^{\frac{m}{2}}+ \cdots +
|\lambda_d|^{\frac{m}{2}} -1 \right)
> \sum_{\frac{m}{k} {\rm prime}, \ k \neq m} |1-(\lambda_1^k+ \cdots
+ \lambda_d^k)|.$$

But if $m$ is large enough then, \begin{equation}|\lambda_1^m +
\cdots + \lambda_d^m|-1
> m\left(|\lambda_1|^{\frac{m}{2}}+ \cdots +
|\lambda_d|^{\frac{m}{2}}+1\right) \label{eqn:m_0} \end{equation}
which is sufficient to give (\ref{eqn:fmgrow}).  For example we
have the inequalities $|1-(\lambda_1^m + \cdots + \lambda_d^m)|>
|\lambda_1|^m -(d-1)|\lambda_2|^m-1$ and
$m\left(1+|\lambda_1|^{\frac{m}{2}}+ \cdots +
|\lambda_d|^{\frac{m}{2}}\right)< m
(d|\lambda_1|^{\frac{m}{2}}+1)$, so whenever,  $$|\lambda_1|^m
-(d-1)|\lambda_2|^m-1> m
\left(d|\lambda_1|^{\frac{m}{2}}+1\right)$$ then (\ref{eqn:m_0})
is satisfied and the proposition is proved for $f \in \MM_b^n$.
When $f \in \MM_{b,k}^n$ for $k<\infty$ and $m \in k\NNN$ then
Remark~\ref{rmk:bfixedL} gives $L(f^m) \leq \#{\rm Fix}(f^m) \leq
2n-1 +L(f^m)$.  Since when $m$ is large, the $2n-1$ term becomes
insignificant in terms of the size of $|\lambda_1|^m$, we see that
we can apply the same proof as above to this case too. \end{proof}

See Examples~\ref{eg:lowgrow} and \ref{eg:m0=3} for applications
of this. Note that there are many examples where the condition
$|\lambda_1|>1$ and $|\lambda_1|> |\lambda_2|$ is not satisfied,
but we still have ${\rm Per}(f) =\NNN$. For example, consider a
map in $\MM_b^2$ which has the action on the first homology of a
matrix with 2's on the diagonal and zeros elsewhere.

Note that in Example~\ref{eg:delaylowgrow} we have a situation
where there are eigenvalues of $f_{\ast 1}$ which are strictly
greater than 1, but ${\rm Per}(f) = 3\NNN$.  So there are limits
to how far we can extend this result.

\begin{rmk}
We would like to estimate $m_0$ for $f \in \MM_b^n$.  From the
above proof, we require that $m_0$ is the infimum of all $m$ \st
$$|\lambda_1|^m
>  d\left[m\left(|\lambda_1|^{\frac{m}{2}}+1\right)+ |\lambda_2|^m\right]+1.$$

\label{rmk:m0}
\end{rmk}

\section{Examples}

\label{sec:examples}

We may apply our results to the following examples.

\begin{eg}
Suppose that $f \in \MM_b^3$ and $f_{\ast 1}$ has matrix
\[
\left( \begin{array}{rrr} 1 & 1 & 1 \\
0 & 0 & 0 \\
1 & 0 & 1
\end{array} \right).
\]
This is a matrix satisfying the conditions of
Proposition~\ref{prop:lowgrow} and so has ${\rm Per}(f) = \NNN$.
Indeed, for $m>1$,

\[
f_{\ast 1}^m = \left( \begin{array}{rrr} 2^{m-1} & 2^{m-2} & 2^{m-1} \\
0 & 0 & 0 \\
2^{m-1} & 2^{m-2} & 2^{m-1}
\end{array} \right)
\]
so we have exponential growth of the trace.  Furthermore, the
eigenvalues of $f_{\ast 1}$ are $0,0$ and $2$. So by
Theorem~\ref{thm:entropy}, the entropy is $\log 2$.
\label{eg:lowgrow}
\end{eg}

The following example has entropy zero.

\begin{eg}
Suppose that $f\in \MM_b^4$ where $f_{\ast 1}$ has action
\[
\left( \begin{array}{rrrr} 1 & 1  & 1 & 1\\
0 & 0 & 0 & 1 \\
0 & 1 & 0 & 0 \\
0 & 0 & 1 & 0
\end{array} \right).
\]
Then we look at the first six iterates of this matrix:

\[
\left( \begin{array}{rrrr} 1 & 1  & 1 & 1\\
0 & 0 & 0 & 1 \\
0 & 1 & 0 & 0 \\
0 & 0 & 1 & 0
\end{array} \right),
\left( \begin{array}{rrrr} 1 & 2  & 2 & 2\\
0 & 0 & 1 & 0 \\
0 & 0 & 0 & 1 \\
0 & 1 & 0 & 0
\end{array} \right),
\left( \begin{array}{rrrr} 1 & 3 & 3 & 3\\
0 & 1 & 0 & 0 \\
0 & 0 & 1 & 0 \\
0 & 0 & 0 & 1
\end{array} \right),
\]

\[
\left( \begin{array}{rrrr} 1 & 4  & 4 & 4\\
0 & 0 & 0 & 1 \\
0 & 1 & 0 & 0 \\
0 & 0 & 1 & 0
\end{array} \right),
\left( \begin{array}{rrrr} 1 & 5  & 5 & 5\\
0 & 0 & 1 & 0 \\
0 & 0 & 0 & 1 \\
0 & 1 & 0 & 0
\end{array} \right),
\left( \begin{array}{rrrr} 1 & 6  & 6 & 6\\
0 & 1 & 0 & 0 \\
0 & 0 & 1 & 0 \\
0 & 0 & 0 & 1
\end{array} \right).
\]

\label{eg:sixcycle} In fact, we see that if $m = 3k$ for $k \geq
1,$ then
\[f_{\ast 1}^m =
\left( \begin{array}{rrrr} 1 & m  & m & m\\
0 & 1 & 0 & 0 \\
0 & 0 & 1 & 0 \\
0 & 0 & 0 & 1
\end{array} \right).
\] and if $3 \nmid m$ then the only non--zero entry on the diagonal
is the top left corner.  Thus $L(f^m)=0$.  If $3|m$ then $L(f^m) =
-3$.  $L(f), L(f^2)=0$.  If $m=3$ then $l(f^3) = 2$.  But when
$m>3$,
$$l(f^m) = \sum_{\stackrel{r|m}{3|\frac{m}{r}}}
\mu(r) L(f^{\frac{m}{r}}).$$ So clearly, if $3\nmid m$ then $l(f)
=0$, and if $3|m$ then $$L(f^m)= -3\sum_{r|\frac{m}{3}}
\mu(r)=0.$$ (For more details of this last calculation, see for
example Section 3 of \cite{eight}.)  Therefore, ${\rm Per}(f) =
\{3\}$.

Also, we calculate that the eigenvalues of this matrix are $1, -1,
e^{\frac{i\pi}{3}}, e^{-\frac{i\pi}{3}}$, and find $l(f^m)$ this
way.  By Theorem~\ref{thm:entropy}, the entropy of this system is
zero (which we would expect since there is no growth of periodic
points).

\end{eg}

We see in the next two examples that a small change to the matrix
in Example~\ref{eg:sixcycle} can alter the entropy and the growth
of periodic points.

\begin{eg}
Now instead consider the matrix in Example~\ref{eg:sixcycle}, but
with any one of the entries $m_{ij}$ for $i, j>1$ which equalled
1, replaced by 2.  Then the eigenvalues of this matrix are $1,
-\left|2^{\frac{1}{3}}\right|, \left|2^{\frac{1}{3}}\right|
e^{\frac{i\pi}{3}}, \left|2^{\frac{1}{3}}\right|
e^{-\frac{i\pi}{3}}$.  By Theorem~\ref{thm:entropy}, the entropy
is $\frac{\log 2}{3}$ and ${\rm Per}(f)= 3\NNN$.  We can see this,
for example, by applying Proposition~\ref{prop:delaylowgrow} to
$f^3$ (since we can compute that the entry in the bottom
right--hand corner of the matrix $f_{\ast 1}^3$ is $2$). We could
also show this by direct calculation.

Note that we cannot apply Proposition~\ref{prop:dominant} here
since $|\lambda_1| = |\lambda_2|$.

\label{eg:delaylowgrow}

\end{eg}

\begin{eg}
Consider $f \in \MM_b^4$ with action
\[f_{\ast 1}=
\left( \begin{array}{rrrr} 1 & 1  & 1 & 1\\
0 & 0 & 0 & 1 \\
0 & 1 & 0 & 0 \\
0 & 0 & 1 & 1
\end{array} \right)
\]
We calculate that
\[f_{\ast 1}^2=
\left( \begin{array}{rrrr} 1 & 2  & 2 & 3\\
0 & 0 & 1 & 1 \\
0 & 0 & 0 & 1 \\
0 & 1 & 1 & 1
\end{array} \right), \
f_{\ast 1}^3=
\left( \begin{array}{rrrr} 1 & 3  & 4 & 6\\
0 & 1 & 1 & 1 \\
0 & 0 & 1 & 1 \\
0 & 1 & 1 & 2
\end{array} \right)
\]
It is easy to see from these matrices that ${\rm Per}(f) = \NNN
\setminus \{2\}$.  Note that we have $m_0=3$ in the statement of
Proposition~\ref{prop:dominant}.  However, note that since the
eigenvalues for $f_{\ast 1}$ are $\lambda_1=1.47, \lambda_2=1,
\lambda_3=0.23+i0.79, \lambda_4=0.23-i0.79$, the calculation in
Remark~\ref{rmk:m0} gives $m_0=10$; which is far from optimal.

\label{eg:m0=3}

\end{eg}

\end{document}